\documentclass[12pt]{article}
\usepackage{local}
\usepackage{geometry}                
\title{Задача Уитни о стержне в поезде, или\\ Много шума из ничего}
\author{А.\,Шень\thanks{Написано по просьбе С.\,М.\,Львовского и В.\,Шувалова, которым автор выражает признательность}}
\date{}

\begin{document}
\maketitle

\begin{flushright}
\emph{Памяти Владимира Андреевича Успенского (1930--2018)}
\end{flushright}

\section{Что пишет Арнольд}

В брошюре 2002 года <<Что такое математика>>~\cite{arnold2002} Владимир Игоревич Арнольд пишет:

\begin{quotation}

Математическая строгость часто оказывается труднопреодолимым препятствием даже
и для хороших математиков.  Следующий пример заимствован из замечательной
классической книги Р.~Куранта и Г.~Роббинса <<Что такое математика?>> (недавно
переизданной на русском языке).

Речь идет о применении топологии.  Пусть на катящейся по горизонтальному
рельсовому пути платформе установлена перпендикулярно рельсам закрепленная
горизонтальная ось, над которой возвышается способный вращаться вокруг этой оси
<<перевернутый маятник>> (стержень).

Утверждается, что \emph{каков бы ни был заданный закон движения платформы \textup(в
течение промежутка времени от нуля до единицы\textup), начальное положение
\textup{<<}маятника\textup{>>} можно выбрать так, что он в конечный момент времени не
будет горизонтален} (Хасслер Уитни).

Авторы доказывают это так. Если исходное положение маятника~--- горизонтально
лежачее, вперед по ходу, то таким оно и останется.  Если же исходное
положение~--- горизонтально лежачее, но назад по ходу, то и это сохранится.

Рассмотрим теперь произвольное начальное положение.  Конечное положение
определяется начальным.  Эта непрерывная функция принимает оба значения
<<вперед>> и <<назад>>.  По теореме топологии она принимает и промежуточные
значения, что и требовалось доказать.

Некоторое время назад мне передали просьбу от проф.~Роббинса (Курант к тому
времени уже умер) постараться исправить это <<ошибочное доказательство>>.  Дело
в том, что никакой непрерывной функции <<конечное состояние при данном
начальном состоянии>> тут сразу не видно:  ее нужно еще точно определить
(с каким-то учетом влияния возможных ударов о платформу), и нужно \emph{доказать}
ее непрерывность.  Я~слышал, что американские математики, пытавшиеся всё это
сделать, написали (неизвестное мне) сочинение с ошибочными промежуточными
утверждениями и доказательствами, так что вопрос о <<маятнике>> и сегодня,
видимо,  остается открытым\footnote
{Дискуссия об этой задаче Х.~Уитни опубликована: B.~E.~Blank.
Book review: What is mathematics? /\!/ Notices of the Amer. Math.
Soc. 2001. Vol.~48, No.~11. P. 1325--1329;
L.~Gillman. Book review: What is mathematics? /\!/ Amer.~Math.~Monthly.
1998. Vol. 105, 
No. 5. P.~485--488; J.~E.~Littlewood. Littlewood's miscellany.
Cambridge Univ. Press, 1986. P.~32--35.    (\emph{Примечание В.\,И.\,Арнольда}. Видимо, в последней ссылке Арнольд имеет в виду~\cite{littlewood1953}.)}.

\end{quotation}

Арнольд возвращается к этому вопросу в 2009 году в другой брошюре, <<Математическое понимание природы. Очерки удивительных физических явлений и их понимания математиками (с рисунками автора)>>~\cite{arnold2009}, где пишет: 

\begin{quotation}

\emph{На горизонтальных рельсах стоит платформа, на которой укреплена перпендикулярно рельсам
горизонтальная ось вращения <<перевёрнутого маятника>>, который может свободно вращаться
в \textup(параллельной рельсам\textup) вертикальной плоскости. Платформа движется по закону
$x=f(t)$ \textup(где $f$~--- гладкая на отрезке $[0, T]$ функция времени\textup).}

\emph{Доказать, что существует такое начальное состояние маятника} $$\alpha(0)=\varphi,\quad \frac{d\alpha}{dt}\,(0)=0,$$
\emph{что он не упадёт на платформу за всё время~$T$ её движения.}

\begin{center}
\includegraphics[scale=.9]{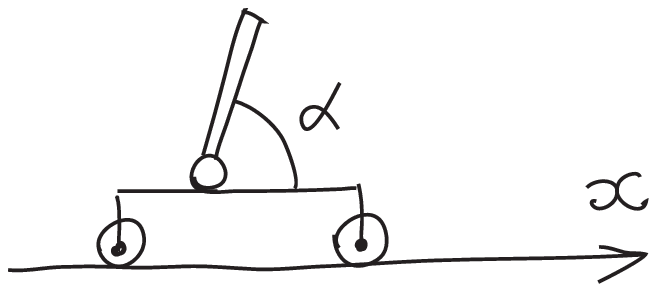}
\end{center}

\textbf{Решение (Куранта)}.  Если $\varphi=0$, то $\alpha(t)=0$ всегда, и~если $\varphi=\pi$,
то $\alpha(t)=\pi$ при всех~$t$.

Из непрерывной зависимости решения (гладкого) дифференциального уравнения от начального условия~$\varphi$
заключаем, по теореме о~промежуточном значении, что между начальными условиями
$\alpha(0)=0$ и~$\alpha(0)=\pi$ существует такое значение $\alpha(0)=\varphi$, что
$\alpha(t)$ при $0\le t \le T$ заключено строго между $0$~и~$\pi$, так что маятник не падает.
$$
\includegraphics[scale=.9]{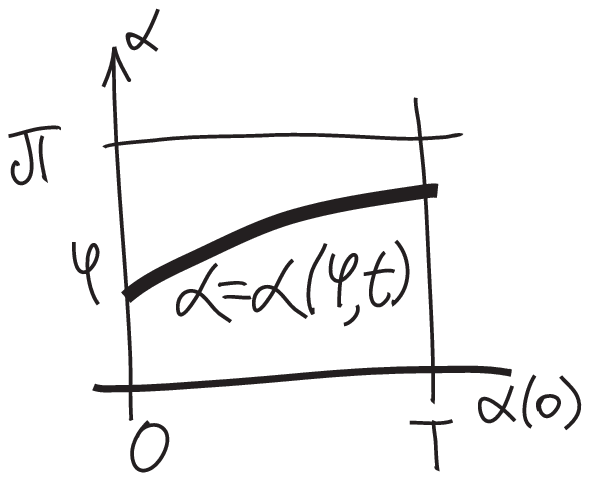}
$$

\textbf{Замечание}.
Многие возражали против этого (ошибочного) доказательства, так как, даже если бы была
определена непрерывная функция $\alpha(\mskip2mu\cdot\mskip2mu, T)$, начального положения~$\varphi$, из её отличия
от нуля и~$\pi$ при начальном условии ${\cdot}=\varphi$ не следует отличие угла~$\alpha$ от~$0$ и~от~$\pi$
во все промежуточные значения моменты времени, $0<t<T$.

Вероятно, можно было бы рассмотреть разумное обобщение рассуждений Куранта, в~котором значения функции
$\alpha (\varphi,\mskip-3mu t)$ определены естественным продолжением после моментов удара о~платформу
(где $\alpha(\varphi, t){=}0$ или~$\pi$). Но в~литературе такое обобщение отсутствует, и~строгое доказательство приведённой гипотезы Куранта неизвестно.

Различные возражения к~различным попыткам обосновать вывод
Куранта опубликовали  Дж.\,Литтлвуд и другие (например,
некоторые из <<контрпримеров>> к~одному из таких обоснований
оказались недействительными, если скорость $df/dt$ меньше
скорости света).

Но я~не видел разумного анализа этой задачи с~учётом ударов.

Описанная теорема была помещена Курантом в~замечательный
элементарный учебник: Курант и~Роббинс <<Что такое
математика?>>, со ссылкой на Уитни.

\end{quotation}

Попытки разобраться в этом деле (в связи с подготовкой переиздания книги Куранта и Роббинса) показали, что тут довольно забавная и запутанная история с участием известных математиков.

\section{Что написано у Куранта и Роббинса}

Первое издание книги <<Что такое математика?>> Куранта и Роббинса~\cite{courant-robbins} вышло на английском языке в 1941 году. В 1947 году вышел его русский перевод под редакцией В.\,Л.\,Гончарова, сделанный по первому и второму английским изданиям (1941 и 1943). 

Действительно, обсуждаемая Арнольдом задача рассматривается в разделе VI.6.2  (2-е издание, с.~352; 3-е издание, с.~347). Авторы пишут (приводим текст по русским изданиям, которые в этой части практически не отличаются друг от друга, разве что небольшой отмеченной ниже стилистической правкой, несколькими запятыми и обозначением угла $\alpha$ на рисунке):

\begin{quotation}
$^*$ \textbf{2. Применение к одной механической проблеме.}
Мы закончим эту главу рассмотрением одной, на первый взгляд трудной, механической проблемы, которая, однако, решается очень просто посредством соображений, связанных с непрерывностью. (Проблема была предложена Г.~У\,и\,т\,н\,и.)

Предположим, что поезд на протяжении некоторого конечного промежутка времени проходит прямолинейный отрезок пути от станции $A$ до станции $B$. Вовсе не предполагается, что движение происходит с постоянной скоростью или с постоянным ускорением. Напротив, поезд может двигаться как угодно: с ускорениями, с замедлениями; не исключены даже мгновенные остановки или частично даже движение в обратном направлении, прежде чем в итоге поезд придет на станцию $B$. Но так или иначе движение поезда на протяжении всего временн\'oго промежутка считается известным заранее; другими словами, считается заданной функция $s=f(t)$, где $s$~--- расстояние поезда от станции $A$, а $t$~--- время, отсчитываемое от момента отправления поезда. К полу одного из вагонов прикреплен на шарнире тяжелый стержень, который без трения может вращаться вокруг оси, параллельной осям вагонов, вперед и назад --- от пола до пола. (Мы допускаем, что, прикоснувшись к полу, он в дальнейшем останется на нем лежать, если ему не случится <<подпрыгнуть>> снова.) [\eng{(If it does touch the floor, we assume that it remains on the floor henceforth; this will be the case if the rod does not bounce.)}\footnote{В английском тексте, как видим, явно указано, что стержень предполагается приклеивающимся к полу при контакте --- а не просто что он будет лежать на полу, <<если ему не случится подпрыгнуть>>.})] 

Вопрос заключается в следующем: \emph{возможно ли в момент отхода поезда поместить стержень в такое начальное положение, т.\,е.~дать ему такой угол наклона, чтобы на протяжении всего пути он не прикоснулся к полу, будучи предоставлен воздействию движения поезда и силы собственной тяжести}?

На первый взгляд может показаться совершенно невероятным, чтобы при наперед определенной схеме движения поезда взаимодействие силы тяжести и сил реакции было способно обеспечить требуемое равновесие стержня при единственном условии --- надлежащем выборе начального положения. Но мы сейчас установим, что такое начальное положение всегда существует.

К счастью, доказательство не подразумевает точного знания законов механики. (Исходя из этих законов, решить задачу было бы чрезвычайно трудно.) Достаточно принять только одно допущение физического содержания: \emph{последующее движение стержня зависит непрерывно от его начального положения}; в частности, если при данном начальном положении стержень во время пути упадет на пол в одну из сторон, то при всяком начальном положении, достаточно мало отличающемся от данного, он не упадет на пол в противоположную сторону.

\begin{center}
\includegraphics[scale=1]{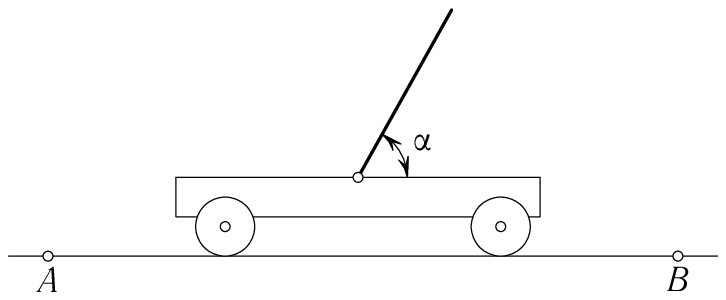}\\
Рис.~175. Проблема Уитни.
\end{center}

Обратим теперь внимание на то, что во всякий момент времени положение стержня характеризуется углом $\alpha$, который он делает с полом [1966 г.; в издании 2001 г. и позже: <<составляет с полом>>]. Углам $\alpha=0^\circ$ и $\alpha=180^\circ$ соответствуют два противоположных горизонтальных (лежачих) положения. Обозначим через $x$ значение угла $\alpha$ в начальном положении стержня. Доказательство нашего утверждения будет косвенное, в соответствии с чисто экзистенциальным характером самой проблемы. Допустим, что \emph{всегда},  т.\,е.~при любом начальном положении стержня, стержень непременно упадет или в одну, или в другую сторону, так что $\alpha$ примет значение или $0^\circ$, или $180^\circ$. Определим тогда функцию $f(x)$ согласно условию:  [1966 г.;  издании 2001 года <<определим тогда функцию $f(x)$ так:>>] $f(x)=+1$ или $-1$, смотря по тому, упадет ли стержень в сторону, соответствующую углу $\alpha=0^\circ$ или углу $\alpha=180^\circ$. Свойства функции $f(x)$ таковы: она задана в интервале $0\le x\le 180^\circ$, непрерывна в нем, и притом $f(0)=+1$, \,$f(180)=-1$. Отсюда, по теореме Больцано, следует, что при каком-то промежуточном значении $x$ ($0^\circ<x<180^\circ$) должно выполняться равенство $f(x)=0$. А это противоречит тому, что функция $f(x)$ может принимать только значения $+1$ и $-1$. Значит, приходится отвергнуть сделанное допущение, согласно которому стержень упадет на пол при каком угодно начальном его положении. 

Совершенно ясно, что приведенное доказательство носит чисто теоретический характер, потому что не дает решительно никаких указаний на то, как определить искомое положение стержня. Вместе с тем, даже если бы такое положение и могло быть вычислено теоретически с абсолютной точностью, практически оно было бы бесполезно вследствие своей неустойчивости. Так, например, в предельном случае, если поезд неподвижен в течении всего <<путешествия>>, решение совершенно очевидно: $x=90^\circ$; но всякий, кто пытался уравновесить иголку в стоячем положении на гладкой горизонтальной поверхности, понимает, насколько это решение практически нереально. Тем не менее с математической точки зрения приведенное доказательство имеет неоспоримый интерес.

\end{quotation} 

\section{Математические подробности}

Проявив занудство (возможно, с элементами бурбакизма), легко восполнить пробелы в изложении Куранта и Роббинса. Но сначала скажем подробнее, в чём состоят эти пробелы (используя базовые сведения об обыкновенных дифференциальных уравнениях).

Движение стержня описывается дифференциальным уравнением (одномерным, второго порядка). Неизвестной функцией является зависимость угла $\alpha$ от времени. Сила, действующая на стержень, складывается из силы инерции (связанной с движением поезда) и силы тяжести. Вместе с силами реакции стержня, позволяющими ему лишь вращаться вокруг оси, они определяют угловое ускорение стержня. Начальное условие говорит, что изначально стержень неподвижен (угловая скорость равна нулю), образуя угол $\alpha$ с горизонтом.

Теорема о непрерывной зависимости решения дифференциального уравнения от начальных условий гарантирует, что значение решения в конце пути поезда (как и в произвольный момент момент времени) является непрерывной функцией от начальных условий. Проблема, отмеченная Арнольдом, состоит в том, что прямо так сослаться на теорему нельзя. Если мы, следуя Куранту и Роббинсу, считаем, что после касания пола стержень <<прилипает>> к нему, то в этот момент дифференциальное уравнение уже не описывает процесс (и вообще ускорение бесконечно, а скорость имеет разрыв, когда стержень с ненулевой скоростью ударяется о пол). Можно убрать пол и рассмотреть движение на всей окружности возможных положений стержня. Тогда уравнение сохраняет силу всегда, но (как правильно отмечает Арнольд) из того, что в конце движения стержень вертикален (или находится в каком-то ином положении выше горизонтали), не следует, что он в процессе движения не опускался ниже горизонтали. Кто знает --- может быть, он даже сделал несколько оборотов, перед тем как вернуться в вертикальное положение.

Как восполнить эти пробелы? Будем, как и предлагалось, рассматривать уравнение на всей окружности (считая, что стержень не ударяется о пол, а может проделать целый оборот). Кроме того, мы можем считать, что задан закон движения платформы при всех $t\ge 0$ (некоторая гладкая функция), и достаточно доказать, что существует начальное условие, при котором решение на всей полуоси $t\ge 0$ остаётся в верхней полуплоскости (в интервале $\pm \pi/2$ от верхнего положения).  Здесь стоит отметить, что стандартные результаты про гладкое уравнение на компактном многообразии гарантируют, что при любом начальном условии решение продолжается на всю полупрямую $t\ge 0$.

Начальное положение в верхней полуплоскости (с нулевой скоростью) будем называть левым (правым), если исходящая из него траектория пересекает горизонтальную полуплоскость первый раз слева (соответственно справа).   Понятие <<первого раза>> определено корректно (если вообще есть пересечение с горизонтальной плоскостью): множество моментов горизонтального положения замкнуто, и среди них есть первый (точная нижняя грань этого множества ему принадлежит). Этот первый момент может быть либо левым, либо правым.  Получаем два множества $L$ и $R$ начальных положений. 

Заметим, что оба множества $L$ и $R$ непусты. В самом деле, функция $f(t)$ в окрестности нуля имеет ограниченную вторую производную, и потому сила инерции, действующая на маятник, ограничена. Следовательно, если начальное положение маятника достаточно близко к горизонтальному, то силы инерции недостаточно, чтобы противодействовать гравитации, и маятник начнёт двигаться вниз и пересечёт горизонталь (с той стороны, где он находился изначально).

Осталось показать, что оба множества $L$ и $R$ открыты, и затем воспользоваться связностью интервала $(0,\pi)$. Почему, скажем, множество $L$ открыто? Пусть $\alpha_0$~--- некоторый элемент $L$,  и $\alpha(t)$~--- соответствующее решение. Пусть $\tau$~--- первый момент горизонтального положения, тогда $\alpha(\tau)=\pi$. На интервале $[0,\tau]$ расстояние от $\alpha(t)$ до другого горизонтального положения ($\alpha=0$) отделено от нуля, поэтому при всех достаточно близких начальных условиях по теореме о непрерывной зависимости (решения как функции от начального условия и момента времени) решение не попадает в правое горизонтальное положение на промежутке $[0,\tau]$. По той же теореме при достаточно близких к $\alpha_0$ начальных условиях решение будет близко к $\pi$ в момент времени $\tau$. Поскольку ускорение тележки ограничено в окрестности отрезка $[0,\tau]$ (и вообще на любом конечном отрезке), то отсюда следует, что при достаточно близких к $\alpha_0$ начальных условиях решение либо уже пересечёт горизонталь слева, либо окажется в таком положении, где это неминуемо. Поэтому эти близкие положения все лежат в $L$, что и требовалось доказать.

\section{Исторические подробности}

\subsection{Пересказ Литлвуда}

Книга Куранта и Роббинса вышла в 1941 году. В 1953 году выходит книга Литлвуда <<Математическая смесь>>~\cite{littlewood1953}, которая в 1962 году была опубликована по-русски в переводе В.\,И.\,Левина. Литлвуд упоминает задачу Уитни, ссылаясь на Куранта и Роббинса. Приведём соответствующий кусок в русском переводе (вполне соответствующем оригиналу):

\begin{quotation}
Стержень шарнирно прикреплён к полу железнодорожного вагона и предоставлен самому себе; тогда существует малая, но отличная от нуля, вероятность того, что он не упадет в течение двух недель: вероятность равна примерно $1:10^{10^5}$. (Поезд не предполагается <<идеальным>>: например, он отправляется со станции Кингс Кросс в 3.15, проходит через туннель, где он 5 минут стоит, а затем, после ряда дальнейших остановок, прибывает в Кембридж в 5.35. Я припоминаю, что мне сообщили, что гений, впервые задавший этот вопрос, не был в состоянии ответить на него.)  [I seem to remember being told that the genius who asked the original question was unable to answer it.] Доказательство имеется в КР [имеется в виду книга Куранта и Роббинса], с.~421--423 [в русском переводе указываются страницы по русскому изданию 1947 года, в английском оригинале Литтлвуд говорит о с.~319].

{\Large$*$} Иной способ доказательства (с разумной свободой интерпретации) состоит в следующем. Рассмотрим начальное положение стержня в относительном покое, составляющее угол $\theta$ с \emph{левым} горизонтальным положением. Пусть $S$ будет множество начальных $\theta$, при которых стержень рано или поздно упадет налево. Основной факт при очень слабых предположениях относительно обстоятельств путешествия и взаимосвязи поезда и стержня (нам нет необходимости детально устанавливать эти предположения) состоит в том, что \emph{множество $S$ открыто}. Пусть $\theta_0$ --- точная верхняя грань углов $\theta$ из $S$. Тогда $\theta_0$ не входит в $S$, и из начального положения $\theta_0$ стержень не упадет налево. Если же он упадет направо, то он сделает это и для всех $\theta$, несколько меньших $\theta_0$; однако это неверно, так как некоторые (в действительности все) из этих $\theta$ принадлежат $S$. Таким образом, из положения $\theta_0$ стержень никогда не упадет. А для некоторого, достаточно малого угла вокруг $\theta_0$ он не упадет в течение двух недель.

Поучительно рассмотреть, почему из этого рассуждения не вытекает, что стержень, соответствующим образом установленный в начальный момент, не будет отклоняться, скажем, более чем на $0{,}5^\circ$ от исходного положения. {\Large$*$}

Это рассуждение моментально убеждает математика. Из всех возможных вариантов избранный мною представляется наиболее подходящим для объяснения любителю. [This argument is instantly convincing to the mathematician. Of possible variants the one chosen seems best suited to be interpreted to the amateur.]

 [Далее Литлвуд предлагает вариант такого более подробного истолкования для <<любителей>>. Для начала он обсуждает понятие точной верхней грани, мы пропускаем этот кусок.]

Предположим, что в начальный момент стержень находится в состоянии относительного покоя и образует угол $\theta$ с \emph{левой} горизонталью; назовем это <<начальным положением $\theta$>>. Рассмотрим множество $S$ начальных положений $\theta$, из которых стержень рано или поздно упадет налево. Если это происходит с ним при каком-нибудь $\theta$, то то же самое будет \emph{иметь место для для всех достаточно близких по обе стороны начальных положений}; любое положение $\theta$, принадлежащее к $S$, находится, следовательно, внутри целого сектора значений $\theta$, принадлежащих к $S$ (на математическом языке: <<$S$ открыто>>). Этот интуитивный факт, на котором все основано, имеет место при весьма общих условиях, уточнять которые нет необходимости. [This intuitive fact, on which everything turns, depends on very slight assumptions, and we need not try ro state them in detail.] 

Пусть $\theta_0$ будет \emph{т.\,в.\,г.} [точной верхней гранью] множества~$S$. Тогда $\theta_0$ не будет элементом $S$, так как в противном случае значения $\theta$, находящиеся справа от $\theta_0$ и достаточно близкие к нему, также принадлежали бы к $S$, и $\theta_0$ не могло бы быть \emph{в.\,г.} [верхней гранью]  множества~$S$. Итак, стержень при начальном положении $\theta_0$ не упадет налево. Но, с другой стороны, если он упадет направо, то он будет падать направо и при \emph{всех} достаточно близких начальных положениях $\theta$ слева от $\theta_0$ (принцип <<открытости>>), так что никакое из этих начальных положений не может быть элементом $S$; но это означает, что $\theta_0$ может быть сдвинуто влево, оставаясь при этом \emph{в.\,г.}~$S$, что противоречит свойству (2) \emph{т.\,в.\,г.} [минимальности]. Следовательно, из исходного положения $\theta_0$ стержень никогда не упадет. А из достаточно малого сектора вокруг $\theta_0$ он не упадет в течении двух недель.

Аналогичный результат верен и при замене обычного шарнира сферическим (позволяющим стержню наклоняться во все стороны). Но в этом случае нам пришлось бы ссылаться уже на значительно более сложную <<теорему о неподвижной точке>> (КР, стр.~335) [в оригинале ссылка на с.~321 английского издания].

\end{quotation}

Ни здесь, ни в другом месте книги Литлвуд не объясняет, каким образом получена оценка длины интервала углов, позволяющего стержню не падать в течение двух недель ($1: 10^{10^5}$). Также не очень понятно, на какого гения намекает Литлвуд
(Уитни?). Можно ещё отметить, что рассуждение (начиная с утверждения об открытости, принятого без доказательства), проведено весьма отчётливо и тут трудно к чему-нибудь придраться даже педанту.

\subsection{Изменение текста Курантом и Роббинсом}

Дополнительную путаницу в истории вопроса создаёт различие в текстах разных изданий Куранта и Роббинса~\cite{courant-robbins}. Мы уже приводили русский перевод, который во всех изданиях следует первому английскому изданию 1941 года. И хотя в английских изданиях текст про задачу Уитни был позже изменён (в издании 1948 года он уже другой), русский текст в последующих переизданиях остался практически неизменным~--- несмотря на заверения А.\,Н.\,Колмогорова в предисловии ко второму русскому изданию, что перевод <<был выправлен и пополнен по последнему английскому (1948) и немецкому (1962) изданиям>>. Тот же самый изменённый текст был вопроизведён в переиздании 1996 года с участием Стьюарта (о комментариях которого смотри ниже). Приводим для справки этот изменённый текст (потому что часть критики Арнольда и других авторов может относиться к нему).

\begin{quotation}

[Формулировка задачи и комментарии к ней оставлены почти без изменений. Начиная с абзаца <<К счастью\ldots>>, изменённый текст выглядит так:] \eng{Paradoxical as this assertion might seem at first sight, it can be proved easily once one concentrates on its essentially topological character.} 

\eng{No detailed knowlegde of the laws of dynamics is needed; only the following simple assumption of the physical nature need be granted: \emph{The motion of the rod depends continuously on its initial position}. Let us characterize the initial position of the rod by the initial angle $x$ which it makes with the floor, and by $y$ the angle which the rod makes with the floor at the end of the journey, when the train reaches the point $B$. If the rod has fallen to the floor we have either $y=0$ or $y=\pi$. For a given initial position $x$ the end position $y$ is, according to our assumption, uniquely determined as a function $y=g(x)$ which is continuous and has the values $y=0$ for $x=0$ and $y=\pi$ for $x=\pi$ (the latter assertion simply expressing that the rod will remain flat on the floor if it starts in this position). Now we recall that $g(x)$, as a continuous function in the interval $0\le x\le \pi$, assumes all the values $y$, e.g. for the value $y=\pi/2$, there exists a specific value of $x$ such that $g(x)=y$; in particular, there exists an initial position for which the end position of the rod at $B$ is perpendicular to the floor.  (Note: in this argument it should not be forgotten that the motion of the train is fixed once for all.)}

\eng{Of course, the reasoning is entirely theoretical. If the journey is of long duration or if the train schedule, expressed by $s=f(t)$, is very erratic, then the range of the initial positions $x$ for which the end position $g(x)$ differs from $0$ or $\pi$ will be exceedingly small, as is known to anyone who has tried to balance a needle upright on a plate for an appreciable time. Still, our reasoning should be of value even to a practical mind inasmuch as it shows how qualitative results in dynamics may be obtained by simple arguments without technical manipulation.}

\end{quotation}

Разница с первоначальным вариантом изложения (см. русский перевод выше) состоит в том, что раньше предлагалось предположить, что стержень всегда падает и рассмотреть функцию с двумя значениями на всех начальных положениях. Непрерывность этой функции означает, что прообразы этих значений открыты, то есть множество тех начальных углов, когда стержень падает влево (или вправо), открыто. Теперь же утверждается безо всяких предположений, что конечное положение стержня (при движении с прилипанием к горизонтали) будет непрерывно зависеть от начальных условий. Это тоже верно, но требует чуть более сложного обоснования (ссылки на теорему о непрерывной зависимости от начального условия тут мало, потому что рассматриваемый вариант движения с прилипанием не соответствует никакому дифференциальному уравнению). Это обоснование также опирается на отмеченный нами выше ключевой факт: попадание в малую зону у горизонтали уже гарантирует скорое попадание на горизонталь (и прилипание), и включает в себя доказательство открытости зон прилипания как часть.

Мне не удалось найти точный момент замены текста (и уж тем более никаких объяснений мотивов такого <<удачного ухудшения>>, если воспользоваться выражением Литлвуда из~\cite{littlewood1953}). Дополнительную путаницу создаёт и то обстоятельство, что перепечатка текста Куранта и Роббинса про задачу Уитни (с заголовком <<\eng{The lever of Mahomet}>> --- видимо, с намёком на одноимённый гроб) в хрестоматии Ньюмана~\cite{courant-robbins-newman} в 1960 году, много позже издания 1948 года с новым текстом, следует оригинальному изданию 1941 года (почему --- не знаю).

\subsection{Разъяснения Бромана}

В 1957 году на 13 Скандинавском математическом конгрессе в Хельсинки Арне Броман сделал доклад, опубликованный в 1958 году~\cite{broman1958}. В этой публикации, со ссылкой на Куранта и Роббинса, формулируется задача Уитни (с требованием доказать, что существует начальное положение, при котором в конце пути стержень будет вертикальным). Приводится решение по Куранту и Роббинсу (в изменённом варианте, с непрерывной зависимостью от начального положения без предположения об обязательном падении).

Далее автор отмечает, что 
<<\eng{However, the physical assumption in the solution of Courant -- Robbins seems a bit hazardous}>>, и делает сноску к этому предложению:  <<\eng{In saying this, we do not intend to criticize the interesting and well-written book by Courant and Robbins. They have included the problem to get an application of the mentioned theorem on continuous functions.}>> После этого автор переходит к обоснованию непрерывности этой самой функции. Сначала он удаляет ограничения на горизонтали:  <<\eng{We now remove for a while the stops for the rod in the horizontal position so that the rod can revolve all around the point $O$. Then the final angle $\beta$, considered as a function $\beta=\beta(\alpha)$ of the initial angle $\alpha$, is obviously defined and continuous for all real~$\alpha$.}>> Таким образом, можно утверждать, что найдётся начальное условие, при котором свободный от горизонтальных ограничений стержень окажется по окончании пути в верхнем положении. Но (как и указывал Арнольд), этого мало --- нужно дополнительно, чтобы стержень всегда оставался в верхней полуплоскости.

Чтобы восполнить этот пробел в рассуждении, Броман замечает, что первое же приближение к горизонтали на достаточно малое расстояние влечёт её пересечение, и рассматривает множество начальных условий, при которых такое пересечение имеет место слева или не имеет места вовсе, но конечное положение находится слева от вертикали. Он устанавливает, что это множество открыто (опираясь на сделанное замечание), и отсюда выводит требуемый результат. 

Наконец, в последнем разделе своей статьи он замечает, что если поезд после приезда в  $B$ навсегда останавливается, то существует такое начальное положение, при котором стержень никогда не коснётся пола (это замечание было уже у Куранта и Роббинса в качестве упражнения).

В качестве примечания 9 Броман приводит немного другое рассуждение (со ссылкой на М.\,Тидемана): в нём дифференциальное уравнение модифицируется так, чтобы стержень, коснувшись пола, продолжал двигаться вниз и не мог сделать полного оборота, но с сохранением достаточной гладкости для применения теоремы о непрерывной зависимости от начальных условий.

\subsection{Критика Постона и её пересказ Стьюартом}

С 1958 года и до 1976 года казалось, что недоразумения с задачей Уитни исчерпаны, и те, кого (вопреки надежде Литлвуда) исходное рассуждение не убеждает, могут обратиться к разъяснениям Бромана. По крайней мере, я не видел никаких публикаций в этот период, авторы которых требовали бы дополнительных разъяснений или выдвигали возражения.

Однако в 1976 году появилась статья Постона в малодоступном издании~\cite{poston1976}. Мне не удалось добыть её текста, но о её содержании можно судить по краткому пересказу самого Постона в его другой статье (см. примечание к описанию статьи~\cite{poston1976}): <<\eng{Certainly there are instances where continuous systems beget discontinuity --- indeed, an interesting case where Courant and Robbins [4] appeal inappropriately to the Intermediate Value Theorem is analysed in [16]}>> (где [16] --- эта самая малодоступная заметка~\cite{poston1976}), и по подробному пересказу в издании Куранта и Роббинса 1996 года~\cite{courant-robbins-stewart}. Это издание вышло с изменениями Стьюарта --- точнее сказать, не с изменениями (поскольку основная часть книги была воспроизведена по предыдущим изданиям), а с дополнениями. Одно из этих дополнений (с.~505 и далее) посвящено задаче Уитни. Стьюарт пишет:

\begin{quotation}
\eng{There is one place where arguably Courant and Robbins made a mistake, although by adding further conditions is it possible to save their argument. Paradoxically, the flaw in their proof is most easily detected if we employ the topological approach to dynamics that their argument was intended to advocate.}

[Далее следует пересказ условия задачи Уитни и решения, опирающегося (после изменений, о которых мы говорили) на непрерывную зависимость конечного положения от начального условия. После этого Стьюарт пишет:]

\eng{The difficulty is that the continuity assumption made in the above discussion is arguably not justified. The problem is not the intricacies of Newton's laws of motion, but those ``absorbing boundary conditions'': if the rod hits the floor, then it stays there. In order to see why the boundary conditions cause trouble, we introduce a topological picture of the possible motions of the system. This approach, known as as \emph{phase portrait}, goes back to Poincar\'e. The idea is to draw a kind of space-time diagram of the motion, not just for a single initial position of the rod, but for many different positions---in principle, all of them. The position of the rod is an angle between $0^\circ$ and $360^\circ$, and we can graph this in the horizontal direction (see Fig.~294). Let time run in the vertical direction. Note that the left and right hand edges of this picture should be identified because $0^\circ=360^\circ$: conceptually, the rectangle is rolled into a cylinder.
\begin{center}
\includegraphics[width=0.7\textwidth]{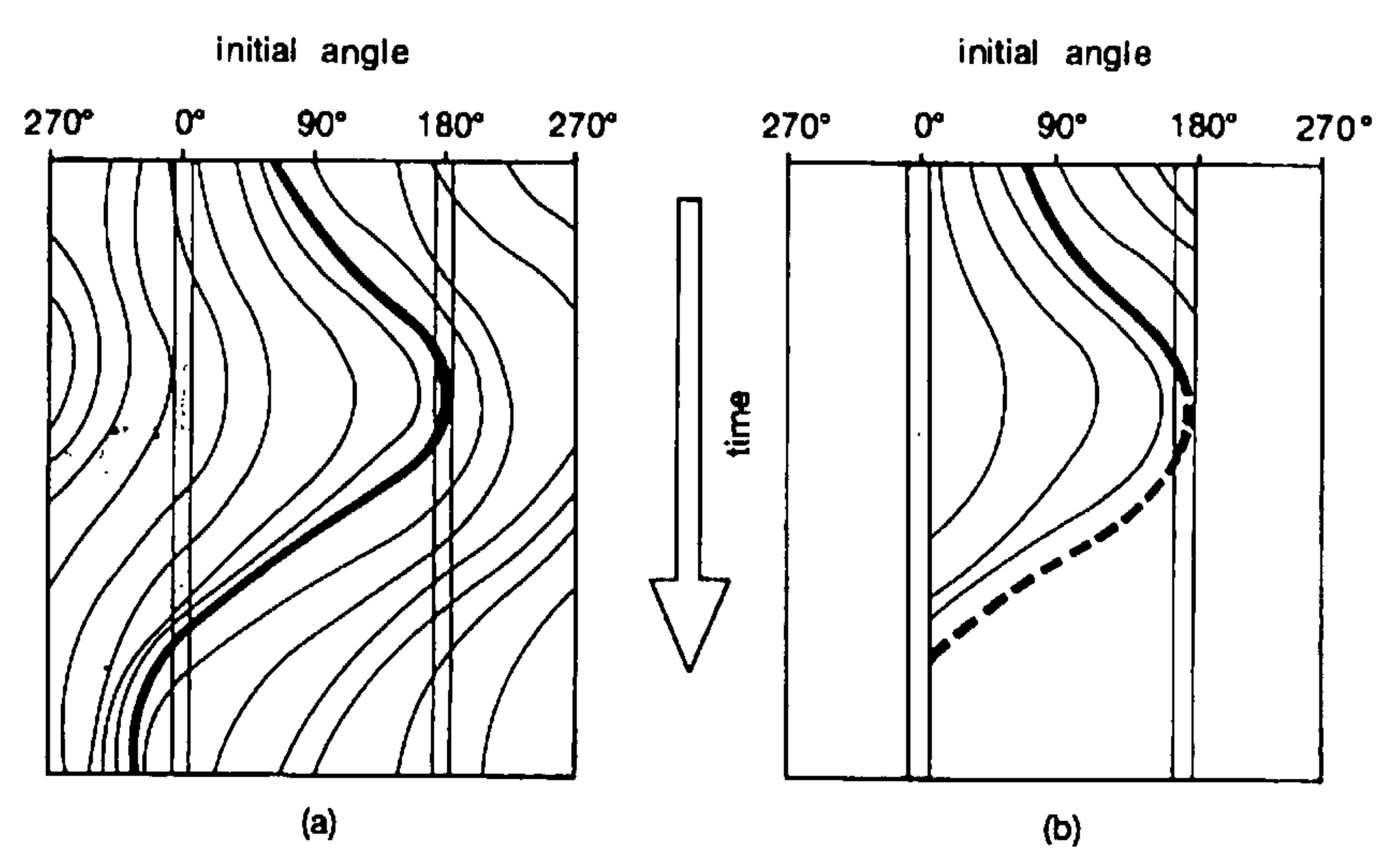}\\
Fig.~294. Possible history of the moving rod for different initial conditions. (a)~Without boundary conditions. (b)~What happens when the boundary conditions are imposed.
\end{center}
}

\eng{Now, the path is space and time of the angle that determines the position of the rod forms some curve that runs up the cylinder --- what Albert Einstein called a ``world-line''. Different initial angles lead to different curves. The laws of dynamics show that these curves vary continuously---provided the boundary conditions are not enforced. Without these conditions the rod is free to turn a full $360^\circ$---there is no floor to prevent it turning all the way round. A possible history is shown in Fig.~294a, and here the final position does depend continuously upon the initial position.}

\eng{However, when the absorbing boundary conditions are put back (Fig.~294b), the final position need \emph{not} depend continuously on the initial one. Curves that just graze the left-hand boundary can swing all the way over to the right. Indeed, in this particular picture \emph{all} initial conditions end up on the floor: contrary to what Courant and Robbins claim, there is no choice that keeps the rod off the floor throughout the motion.}

\eng{This error in Courant and Robbins' reasoning was first pointed out by Tim Poston in 1976, but it is still not widely known. The continuity assumption can be resuscitated by imposing extra con[s]traints on the motion, for example a perfectly level track, no springs on the train, and so forth. But it seems more instructive, as an exercise in the application of topology to dynamics, to understand why the absorbing boundary conditions destroy continuity.}

\end{quotation}

С математической точки зрения Стьюарт пишет правильно: действительно, добавление поглощающих границ нарушает непрерывность, и действительно в предположении о том, что поезд движется горизонтально и не имеет рессор, непрерывность может быть обоснована.  Примерно такой же (но чуть менее подробный) текст можно найти в книге Стьюарта~\cite{stewart1989}, где он, однако, не упоминает о том, что в предположениях Куранта и Роббинса непрерывность несложно доказать, а говорит лишь, что <<\eng{the continuity assumption is not justified}>> (и приводит тот же рисунок траекторий без граничных условий и с ними).

Однако и то верно, что у Куранта и Роббинса эти дополнительные предположения были сделаны, а утверждение о непрерывности (для этого, а не более общего случая) было сформулировано без обоснования (и это было явно указано). То, что для негоризонтального пола (поглощающих границ под произвольным углом) утверждение будет неверно, отмечал ещё Литлвуд (абзац <<Поучительно рассмотреть>>  в приведённой выше цитате). В такой ситуации говорить об <<ошибке в рассуждении Куранта и Роббинса>>, как это делают Стьюарт и Постон (см. его цитату в~\cite{gillman1998}), кажется несколько чрезмерным.

Похоже, что впоследствии и сам Стьюарт пришёл к такому заключению. По крайней мере, такое впечатление создаётся при чтении отчёта о его лекции в 2009 году~\cite{stewart2009}:

\begin{quote}
\eng{Professor Stewart’s lecture, entitled \emph{The Strange Case of the Courant -- Robbins Train}, centered around a problem that he said had been bothering him for 40 years. It concerns a train and how its motion, and the influence of gravity, would affect a rod pivoted on the floor. In their 1941 book \emph{What is Mathematics?} Richard Courant and Herbert Robbins asked if it were possible to place the rod in such a position that, if it is released when the train starts, it will not fall to the floor during the entire journey. Although this appears to be a complex problem in mechanics, they asserted that the answer is a simple yes and avoided the complexities of mechanics using a simple application of the concept of continuity. With his trademark geniality and witty presentation skills, Professor Stewart carefully guided the audience through his concern: does the concept of continuity actually apply in this situation? By considering boundary conditions and floors of railway carriages which slanted towards the middle, he showed that in fact Courant and Robbins were correct to assume continuity in the particular case where the carriage has a flat floor although the proof of continuity actually does require a careful analysis of the underlying mechanics. But this is not always the case for other shaped floors.}
\end{quote}

\subsection{Рецензии на новое издание}

Выход нового издания книги Куранта и Роббинса с дополнениями Стьюарта был замечен  рецензентами. В 1998 году в \emph{Americal Mathematical Monthly} была опубликована рецензия Леонарда Гиллмана~\cite{gillman1998}. В ней приводится формулировка задачи и излагаются возражения Стьюарта со ссылкой на Постона. (Гиллману, в отличие от меня, удалось получить текст статьи Постона~\cite{poston1976} --- ему его прислал Мартин Гарднер!) Реакция Постона Гиллману тоже кажется чрезмерной (он отмечает, что <<\eng{Poston's term is `wrong'}>>, p.~488); далее он, со ссылкой на Радина (Charles Radin, частное сообщение) намечает отсутствующее у Куранта и Роббинса рассуждение. Радин также использует ключевое соображение о том, что приблизившись к полу достаточно близко, стержень непременно столкнётся с ним (видимо, в тексте опечатка в знаке неравенства для~$t'$):

\begin{quote}
\eng{First we prove: For any given physical setup and itinerary $x(t)$, there is an $\epsilon>0$ such that, if $H(t)<\epsilon$ and $V(t)$ is either zero or directed toward the floor, then $H(t')\le H(t)$ for $t'<t<T$. This follows almost immediately from noticing that when the rod is almost at the floor, the angular acceleration of the rod produced by the bounded horizontal acceleration of the train must be smaller in magnitude that the vertical acceleration caused by gravity; the net acceleration must then be directed toward the floor, and if the angular verlocity is also directed toward the floor the conclusion follows.}
\end{quote}

Далее он объясняет (примерно так же, как и Броман), почему из этого соображения вытекает непрерывность.

В другой рецензии (Брайн Бланк~\cite{blank2001}, \emph{Notices of the AMS}) рецензент с сожалением отмечает, что изменения в издании 1996 года по сравнению с предыдущими минимальны: 

\begin{quote}
\eng{All the old figures are there. Even the original typesetting has been preserved --- quite a welcome change from the homogenized TEX of the present. Closer inspection, however, reveals something less welcome: the new edition is by and large a photographic reproduction. I could not find even one change in it. Stewart rationalizes this in his preface: ``not a single word or symbol had to be deleted from this new edition.'' Really? After half a century could we not have had a correction to the name of Mr. Arch medes (p. 400)?}
\end{quote}

Рецензент также упоминает критику Стьюарта и объяснения Гиллмана:

\begin{quote}
\eng{Stewart relates a challenge to the basic continuity assumption that underlies the solution given by Courant and Robbins. Ironically, that challenge itself relies on the assumption that the train can move in a way that is physically impossible, as Gillman has pointed out in his review [5].}
\end{quote}

(Здесь [5] --- это рецензия Гиллмана~\cite{gillman1998}.)

\subsection{После Арнольда}

После публикаций книг Арнольда (который, кстати, даёт ссылку на рецензии Гиллмана и Бланка, но никак их не комментирует), интерес к задаче Уитни снова вырос, и появилось несколько сравнительно недавних публикаций~\cite{polekhin2014,polekhin2014a,zubelevich2015,bolotin-kozlov2015,srzednicki2017}.

И.\,Ю.\,Полехин в статье 2014 года~\cite{polekhin2014}, приведя формулировку задачи Уитни из русского издания Куранта и Роббинса, пишет:
\begin{quote}
Положительный ответ на этот вопрос [о начальном положении, когда стержень не коснётся пола] дается самими авторами. Тем не менее, предложенное доказательство использует предположение о том, что движение стержня зависит непрерывным образом от начальных условий. Данное предположение кажется естественным, и подобное действительно верно для широкого класса механических систем, но в случае, когда считается, что, коснувшись пола, стержень вечно остается горизонтальным, непрерывная зависимость от начальных данных становится менее очевидной. В частности, эта неполнота доказательства была отмечена Арнольдом при обсуждении данной задачи [2]. В настоящей работе приводится строгое обоснование существования решения без падений (в том числе для случая решения, определенного на бесконечном полуинтервале времени), основанное на применении топологического метода Важевского.
\end{quote}
Здесь [2] --- это книга Арнольда~\cite{arnold2002}. Само рассуждение, насколько я могу судить (работ Важевского я не читал, и в изложении Полехина много формул), следует схеме Бромана и Гиллмана. Помимо этого, доказывается, что для случая периодического движения поезда существует периодическое решение для стержня, при котором он никогда не касается пола.

В английском варианте~\cite{polekhin2014a} Полехин описывает ситуацию так:
\begin{quote}
\eng{This shortcoming of the original proof} [из книги Куранта и Роббинса]  \eng{was mentioned and briefly commented by Arnold in his book [Arn] which also includes short overview of the articles related to the matter. Yet detailed consideration of the issue of whether original prove} [опечатка: вероятно, должно быть `proof'] \eng{is full and correct or not is beyond the objective of the paper. In the first section, application of the topological Wa\.zewski principle to the above problem is considered. This method allows us to prove existence of solutions without falling (even on an infinite time interval), including solutions with zero initial velocity of inverted pendulum.}
\end{quote}
Здесь [Arn] --- та же книга Арнольда~\cite{arnold2002}. Полехин упоминает статьи, цитируемые Арнольдом, включая рецензии Гиллмана и Бланка, но никак не комментирует связи своего рассуждения с приводимым Гиллманом рассуждением Радина, и осторожно пишет про <<\eng{issue of whether original proof is full and correct or not}>>. Статья Бромана~\cite{broman1958} также не упоминается.

Есть ещё одна публикация примерно того же времени, посвящённая задаче Уитни и её обобщениям~\cite{zubelevich2015} --- её автор Олег Зубелевич, изложив условие задачи Уитни, пишет:
\begin{quote}

\eng{The authors} [Курант и Роббинс] \eng{gave positive answer to this question. Their argument was informal. V. Arnold in [1] considered this problem as open. The complete solution to the problem has been given by I. Polekhin in his Ph.D. thesis (unpublished), see also [3]. He solved the problem by direct application of results from [4]. In this article we prove simple and general theorem which implies particularly that there are continuum never-falling solutions to Whitney’s problem, we also believe that this theorem describes many other such a type effects.}
\end{quote}
Здесь [1] --- та же книга Арнольда 2002 года, [3] --- публикация Полехина~\cite{polekhin2014a}, а [4] --- статья~\cite{srzednicki2005} (где непосредственно задача Уитни не упоминается). Один из авторов последней статьи, Роман Сшедницки, в 2017 году опубликовал заметку~\cite{srzednicki2017}, в которой рассказывает об истории задачи Уитни (и из которой я узнал о существовании работы Бромана~\cite{broman1958}, в частности). Он пишет:

\begin{quote}
\eng{In the book ``What is Mathematics?'' Richard Courant and Herbert Robbins presented a solution of a Whitney’s problem of an inverted pendulum on a railway carriage moving on a straight line. Since the appearance of the book in 1941 the solution was contested by several distinguished mathematicians. The first formal proof based on the idea of Courant and Robbins was published by Ivan Polekhin in 2014. Polekhin also proved a theorem on the existence of a periodic solution of the problem provided the movement of the carriage on the line is periodic. In the present paper we slightly improve the Polekhin’s theorem by lowering the regularity class of the motion and we prove a theorem on the existence of a periodic solution if the carriage moves periodically on the plane.}
\end{quote}

Однако, говоря о статье Бромана~\cite{broman1958}, Сшедницки довольно загадочно пишет <<\eng{He presented a comprehensive argument supporting the Courant and Robbins’ solution, although his proof lacks of formal rigor at some details}>> --- не указывая, однако, что же это за такие детали. Упомянув возражения Арнольда в~\cite{arnold2009} (<<\eng{At that time Arnold had still objections towards the correctness of the solution from [CR]}>>), он пишет:

\begin{quote}
\eng{Finally, in 2014, 73 years from the announcement of the Whitney’s problem, in the paper [P1] Ivan Polekhin provided a short rigorous proof of the Courant and Robbins’s solution based on the Wa\.zewski retract theorem. Other proofs were published in [BK] and [Zu].}
\end{quote}

Здесь [P1] --- это статьи Полехина~\cite{polekhin2014,polekhin2014a}, [BK] --- это статья Болотина и Козлова~\cite{bolotin-kozlov2015},  а [Zu] --- статья Зубелевича~\cite{zubelevich2015}. В статье~\cite{bolotin-kozlov2015} (после формулировки некоторых общих результатов) упоминается задача Уитни, ситуация с которой описывается так:

\begin{quote}
Утверждение о существовании движений перевернутого маятника без падений сформулировано Уитни и обсуждается в книге [2]. Доказательство основано на несуществовании ретракции отрезка на его границу. Критические замечания В. И. Арнольда [5] стимулировали появление полных доказательств теоремы Уитни (с указанием условий ее справедливости) [3], [4]. 
\end{quote}

Здесь [2] --- это книга Куранта и Роббинса, несуществование ретракции --- теорема о промежуточных значениях, [5] --- книга~\cite{arnold2002}, а [3] и [4] --- соответственно~\cite{polekhin2014} и~\cite{zubelevich2015}. Не вполне ясно, считают ли авторы~\cite{bolotin-kozlov2015} доказательства Бромана из~\cite{broman1958} или Радина из~\cite{gillman1998} неполными и если да, то в чём их неполнота.

В заключение скажу ещё о своём разговоре с Сабиром Меджидовичем Гусейн-Заде на одной из летних школ в Дубне, где я услышал о возражениях Арнольда по поводу рассуждения в Куранте и Роббинсе. Арнольд был тогда ещё жив, но подойти к нему с этим я (как и многие другие) боялся. Вместо этого я настойчиво просил Гусейн-Заде послушать, как я представлял себе доказательство Куранта -- Роббинса (примерно в том же духе, как написано выше), и сказать мне отчётливо, правильно ли это рассуждение. Вежливый Сабир Меджидович меня выслушал и сказал, что ошибки в рассуждении не видит (но, кажется, с Арнольдом это после не обсуждал, как и я).

\end{document}